\documentclass{article}
\usepackage{amsfonts,amssymb}
\usepackage{amsmath}

\begin{document}

\title{Lebesgue constants on compact manifolds}
\author{A. Kushpel \\
Department of Mathematics,\\
University of Leicester}
\date{\today}
\maketitle

\begin{center}
\textbf{Abstract}
\end{center}

Sharp asymptotic for norms of Fourier projections on compact homogeneous
manifolds $\mathbb{M}^{d}$ (for example, the real spheres $%
\mathbb{S}^{d}$, the real, complex and
quaternionic projective spaces $\mathrm{P}^{d}(\mathbb{R})$, $\mathrm{P}^{d}(%
\mathbb{C})$, $\mathrm{P}^{d}(\mathbb{H})$ and the Cayley elliptic plain $%
\mathrm{P}^{16}(\mathrm{Cay})$) are established. These results extend sharp
asymptotic estimates found by Fejer \cite{fejer} in the case of $\mathbb{S}%
^{1}$ in 1910 and then by Gronwall \cite{gronwall} in 1914 in the case of $%
\mathbb{S}^{2}$. As an application of these results we give solution of the
problem of Kolmogorov on sharp asymptotic for the rate of convergence of
Fourier sums on a wide range of sets of multiplier operators.\\



\section{Introduction}

Let $\mathbb{M}^{d}$ be a compact globally symmetric space of rank 1
(two-point homogeneous space), $\nu $ its normalized volume element, $\Delta
$ its Laplace-Beltrami operator. It is well-known that the eigenvalues $%
\theta _{k}$, $k\geq 0$ of $\Delta $ are discrete, nonnegative and form an
increasing sequence $0\leq \theta _{0}\leq \theta _{1}\leq \cdots \leq
\theta _{k}\leq \cdots $ with $+\infty $ the only accumulation point.
Corresponding eigenspaces $\mathrm{H}_{k}$, $k\geq 0$ are finite
dimensional, $d_{k}=\mathrm{dim}\mathrm{H}_{k}<\infty $, $k\geq 0$,
orthogonal and $L_{2}(\mathbb{M}^{d},\nu )=\oplus _{k=0}^{\infty }\mathrm{H}%
_{k}$. Let $\{Y_{j}^{k}\}_{j=1}^{d_{k}}$ be an orthonormal basis of $\mathrm{%
H}_{k}$. Assume that $\phi $ is a continuous function $\phi \in C(\mathbb{M}%
^{d})$ with the formal Fourier expansion
\[
\phi \sim \sum_{k=0}^{\infty }\sum_{j=1}^{d_{k}}c_{k,j}(\phi
)Y_{j}^{k},\,\,\,c_{k,j}(\phi )=\int_{\mathbb{M}^{d}}\phi \overline{Y_{j}^{k}%
}d\nu .
\]%
Consider the sequence of Fourier sums
\[
S_{n}(\phi )=\sum_{k=0}^{n}\sum_{j=1}^{d_{k}}c_{k,j}(\phi )Y_{j}^{k}.
\]%
The main aim of this article is to establish sharp asymptotic for the
sequence $\Vert S_{n}\Vert _{C(\mathbb{M}^{d})\rightarrow C(\mathbb{M}^{d})}$
as $n\rightarrow \infty $. As a consequence we give the solution of the
problem of Kolmogorov on sharp asymptotic for the rate of convergence of
Fourier sums on sets generated by pseudo-differential operators (multiplier
operators) on compact manifolds.

Observe that this set of problems is closely related to the problem of
uniform convergence of Fourier series on $\mathbb{M}^{d}$. Indeed, let
\[
E_{n}(\phi) = \inf \left\{\|\phi-t_{n}\|_{C(\mathbb{M}^{d})} \left|\right.
t_{n} \in \mathcal{T}_{n}\right\}
\]
be the best approximation of a function $\phi \in C(\mathbb{M}^{d})$ by the
subspace $\mathcal{T}_{n}$ of polynomials of order $\leq n$, $\mathcal{T}%
_{n} = \bigoplus_{k=0}^{n} \mathrm{H}_{k}$. Then, by the Lebesgue inequality
\cite{lebesgue} we get
\[
\|\phi - S_{n}(\phi)\|_{C(\mathbb{M}^{d})} \leq \left(1+\|S_{n}\|_{C(\mathbb{%
M}^{d}) \rightarrow C(\mathbb{M}^{d})}\right)E_{n}(\phi),
\]
where $\|S_{n}\|_{C(\mathbb{M}^{d}) \rightarrow C(\mathbb{M}^{d})} = \sup
\{\|S_{n}(\phi)\|_{C(\mathbb{M}^{d})} |\, \|\phi\|_{C(\mathbb{M}^{d})} \leq
1\}$. It means that $S_{n}(\phi, x)$ converges uniformly to $\phi$ if
\[
\lim_{n \rightarrow \infty}\,\,E_{n}(\phi)\|S_{n}\|_{C(\mathbb{M}^{d})
\rightarrow C(\mathbb{M}^{d})}\,=\,0.
\]
In the case of the circle, $\mathbb{S}^{1}$, the following result has been
found by Fejer in 1910 \cite{fejer},
\[
\|S_{n}\|_{C(\mathbb{S}^{1}) \rightarrow C(\mathbb{S}^{1})} = \frac{1}{\pi}%
\,\int_{-\pi}^{\pi} |D_{n}(t)| dt \,= \frac{4}{\pi^{2}}\,\, \mathrm{ln} \,n
+ O(1).
\]
where $D_{n}(t) = 1/2 + \sum_{k=1}^{n} \cos\,\,kt$ is the Dirichlet kernel.
In the case of\,\, $\mathbb{S}^{2}$, the two-dimensional unit sphere in $%
\mathbb{R}^{3}$, the estimates of $\|S_{n}\|_{C(\mathbb{S}^{2}) \rightarrow
C(\mathbb{S}^{2})}$ as $n \rightarrow \infty$, have been established by
Gronwall \cite{gronwall}. Namely, it was shown that
\[
\|S_{n}\|_{C(\mathbb{S}^{2}) \rightarrow C(\mathbb{S}^{2})}= n^{1/2}\,\,%
\frac{2}{\pi^{3/2}}\int_{0}^{\pi} \sqrt{\mathrm{cot}\left(\frac{\eta}{2}%
\right)}\,\,d\eta + O(1),\,
\]
\[
=\, n^{1/2} 2^{3/2} \pi^{-1/2}+O(1).
\]

\section{Harmonic Analysis on Compact Manifolds}

We shall be mostly interested here in compact globally symmetric spaces of
rank 1 (two-point homogeneous spaces) and the complex sphere $\mathbb{S}_{%
\mathbb{C}}^{d}$. Such manifolds of dimension $d$ will be denoted by $%
\mathbb{M}^{d}$. In particular, each $\mathbb{M}^{d}$ can be considered as
the orbit space of some compact subgroup $\mathcal{H}$ of the orthogonal
group $\mathcal{G}$, that is $\mathbb{M}^{d}=\mathcal{G}/\mathcal{H}$. Let $%
\pi :\mathcal{G}\rightarrow \mathcal{G}/\mathcal{H}$ be the natural mapping
and $\mathbf{e}$ be the identity of $\mathcal{G}$. The point $\mathbf{o}=\pi
(\mathbf{e})$ which is invariant under all motions of $\mathcal{H}$ is
called the pole (or the north pole) of $\mathbb{M}^{d}$. On any such
manifold there is an invariant Riemannian metric $d(\cdot ,\cdot )$, and an
invariant Haar measure $d\nu $. Two-point homogeneous spaces admit
essentially only one invariant second order differential operator, the
Laplace-Beltrami operator $\Delta $. A function $Z(\cdot ):\mathbb{M}%
^{d}\rightarrow \mathbb{R}$ is called zonal if $Z(h^{-1}\cdot )=Z(\cdot )$
for any $h\in \mathcal{H}$. A complete classification of the two-point
homogeneous spaces was given by Wang \cite{wang}. For information on this
classification see, e.g., Cartan \cite{car}, Gangolli \cite{gang}, and
Helgason \cite{hel1,helgason}. The geometry of these spaces is in many
respects similar. All geodesics in a given one of these spaces are closed
and have the same length $2L$. Here $L$ is the diameter of $\mathcal{G}/%
\mathcal{H}$, i.e., the maximum distance between any two points. A function
on $\mathcal{G}/\mathcal{H}$ is invariant under the left action of $\mathcal{%
H}$ on $\mathcal{G}/\mathcal{H}$ if and only if it depends only on the
distance of its argument from $\mathbf{o}=\mathbf{e}\mathcal{H}$. Since the
distance of any point of $\mathcal{G}/\mathcal{H}$ from $\mathbf{e}\mathcal{H%
}$ is at most $L$, it follows that a $\mathcal{H}$-spherical function $Z$ on
$\mathcal{G}/\mathcal{H}$ can be identified with a function $\tilde{Z}$ on $%
[0,L]$. Let $\theta $ be the distance of a point from $\mathbf{e}\mathcal{H}$%
. We may choose a geodesic polar coordinate system $(\theta ,\mathbf{u})$
where $\mathbf{u}$ is an angular parameter. In this coordinate system the
radial part $\Delta _{\theta }$ of the Laplace-Beltrami operator $\Delta $
has the expression
\[
\Delta _{\theta }=(A(\theta ))^{-1}\frac{d}{d\theta }\left( A(\theta )\frac{d%
}{d\theta }\right) ,
\]%
where $A(\theta )$ is the area of the sphere of radius $\theta $ in $%
\mathcal{G}/\mathcal{H}$ which can be computed in terms of the structure of
the Lie algebras of $\mathcal{G}$ and $\mathcal{H}$ (see Helgason \cite[p.251%
]{helgason}, \cite[p.168]{hel1} for the details). It can be shown that
\begin{equation}
A(\theta )=\omega _{\sigma +\rho +1}\lambda ^{-\sigma }(2\lambda )^{-\rho
}(\sin \lambda \theta )^{\sigma }(\sin 2\lambda \theta )^{\rho },
\label{a theta}
\end{equation}%
where $\omega _{d}$ is the area of the unit sphere in $\mathbb{R}^{d}$ and

\[
\mathbb{S}^{d}:{\ }\sigma =0,\rho =d-1,\lambda =\pi /2L,d=1,2,3,\ldots ;
\]%
\[
\mathrm{P}^{d}(\mathbb{R}):{\ }\sigma =0,\rho =d-1,\lambda =\pi
/4L,d=2,3,4,\ldots ;
\]%
\begin{equation}  \label{pipa}
\mathrm{P}^{d}(\mathbb{C}):{\ }\sigma =d-2,\rho =1,\lambda =\pi
/2L,d=4,6,8,\ldots ;
\end{equation}%
\[
\mathrm{P}^{d}(\mathbb{H}):{\ }\sigma =d-4,\rho =3,\lambda =\pi
/2L,d=8,12,\ldots ;
\]%
\[
\mathrm{P}^{16}(\mathrm{Cay}):{\ }\sigma =8,\rho =7,\lambda =\pi /2L.
\]
Applying (\ref{a theta}) and (\ref{pipa}) we can write the operator $\Delta
_{\theta }$ (up to some numerical constant) in the form
\[
\Delta _{\theta }=\frac{1}{(\sin \lambda \theta )^{\sigma }(\sin 2\lambda
\theta )^{\rho }}\frac{d}{d\theta }(\sin \lambda \theta )^{\sigma }(\sin
2\lambda \theta )^{\rho }\frac{d}{d\theta }.
\]%
Using a simple change of variables $t=\cos 2\lambda \theta $, this operator
takes the form (up to a positive multiple),
\begin{equation}
\Delta _{t}=(1-t)^{-\alpha }(1+t)^{-\beta }\frac{d}{dt}(1-t)^{1+\alpha
}(1+t)^{1+\beta }\frac{d}{dt},  \label{diana}
\end{equation}%
where
\begin{equation}
\alpha =\frac{\sigma +\rho -1}{2},\hspace{1cm}\beta =\frac{\rho -1}{2}.
\label{pipo}
\end{equation}%
For all manifolds considered here $\alpha =(d-2)/2.$ We will need the
following statement Szeg\"{o} \cite[p.60]{sze}:

\textbf{Proposition 1.} \emph{The Jacobi polynomials $y=P_{k}^{(\alpha,
\beta)}$ satisfy the following linear homogeneous differential equation of
the second order:
\[
(1-t^{2})y^{^{\prime \prime }}+(\beta-\alpha-(\alpha+\beta+2)t)y^{^{\prime
}}+ k(k+\alpha+\beta+1)y=0,
\]
or
\[
((1-t)^{\alpha+1}(1-t)^{\beta+1}y^{^{\prime }})^{^{\prime }}+
k(k+\alpha+\beta+1)(1-t)^{\alpha}(1+t)^{\beta}y = 0.
\]
} It follows from the above proposition that the eigenfunctions of the
operator $\Delta_{t}$ defined in (\ref{diana}) are well-known Jacobi
polynomials $P_{k}^{(\alpha, \beta)}$ and the corresponding eigenvalues are $%
\theta_{k} = -k(k+\alpha+\beta+1)$. In this way zonal $\mathcal{H}$%
-invariant functions on $\mathbb{M}^{d}=\mathcal{G}/\mathcal{H}$ can be
easily identified in each of the five cases above since the elementary zonal
functions are eigenfunctions of the Laplace-Beltrami operator. We shall call
them $Z_k$, $k \in \mathbb{N} \cup \{0\}$, with $Z_0 \equiv 1$. Let $\tilde{Z%
}_k$ be the corresponding functions induced on $[0, L]$ by $Z_k$. Then
\begin{equation}  \label{z}
\tilde{Z}_k(\theta) = C_{k}(\mathbb{M}^{d}) P_{k}^{(\alpha, \beta)}(\cos
2\lambda \theta),\,\,\,k \in \mathbb{N} \cup \{0\},
\end{equation}
where \,$\alpha$\, and \,$\beta$\, have been specified above. If $\mathbb{M}%
^{d} = \mathrm{P}^{d}(\mathbb{R})$, then only the polynomials of even degree
appear because, due to the identification of antipodal points on $\mathbb{S}%
^{d}$, only the even order polynomials $P_{k}^{(\alpha, \alpha)}$, $k=2m$, $%
m \in \mathbb{N} \cup \{0\}$, can be lifted to be functions on $\mathrm{P}%
^{d}(\mathbb{R})$.

In the case of $\mathbb{S}^{d}$ we have $\sigma=0$, $\rho=d-1$, so that, $%
\alpha=\beta=(d-2)/2$ and the polynomials $P_{k}^{(\alpha, \beta)}$ reduce
to $P_{k}^{((d-2)/2, (d-2)/2)}$ which is a multiple of the Gegenbauer
polynomial $P_{k}^{(d-1)/2}$. A detailed treatment of the Jacobi polynomials
can be found in Szeg\"o \cite{sze}. We remark that the Jacobi polynomials $%
P_{k}^{(\alpha, \beta)}(t)$, $\alpha > -1, {\ } \beta > -1$ are orthogonal
with respect to $\omega^{\alpha, \beta}(t)=c^{-1}(1-t)^{\alpha}(1+t)^{\beta}$
on $(-1, 1)$. The above constant $c$ can be found using the normalization
condition $\int_{\mathbb{M}^{d}} d\nu = 1$ for the invariant measure $d\nu$
on $\mathbb{M}^{d}$ and a well-known formula for the Euler integral of the
first kind
\begin{equation}  \label{euler}
\mathrm{B}(p,q) = \int_{0}^{1} \xi^{p-1} (1-\xi)^{q-1}d\xi = \frac{%
\Gamma(p)\Gamma(q)}{\Gamma(p+q)},\,\,\,p>0,\,\,q>0.
\end{equation}
Applying (\ref{euler}) and a simple change of variables we get
\[
1= \int_{\mathbb{M}^{d}} d\nu = \int_{-1}^{1} \omega^{\alpha, \beta}(t)dt =
c^{-1} \int_{-1}^{1} (1-t)^{\alpha}(1+t)^{\beta} dt,
\]
so that,
\begin{equation}  \label{c}
c = \int_{-1}^{1} (1-t)^{\alpha}(1+t)^{\beta} dt = 2^{\alpha + \beta + 1}
\frac{\Gamma(\alpha + 1) \Gamma( \beta + 1)}{\Gamma(\alpha + \beta + 2)}.
\end{equation}
We normalize the Jacobi polynomials as follows:
\[
P_{k}^{(\alpha, \beta)}(1) = \frac{\Gamma(k+ \alpha +1)}{\Gamma(\alpha +
1)\Gamma(k+ 1)}.
\]
This way of normalization is coming from the definition of Jacobi
polynomials using the generating function Szeg\"o \cite[p.69]{sze}.

Let $L_{p}(\mathbb{M}^{d})$ be the set of functions of finite norm given by
\[
\parallel \varphi \parallel _{p}\,\,=\,\,\parallel \varphi \parallel _{L_{p}(%
\mathbb{M}^{d})}=\left\{
\begin{array}{c}
\left( \int_{\mathbb{M}^{d}}\left\vert \varphi (x)\right\vert ^{p}d\nu
(x)\right) ^{1/p},1\leq p<\infty , \\
\mathrm{ess}\sup \left\vert \varphi \right\vert ,p=\infty .%
\end{array}%
\right.
\]%
Further, let $U_{p}=\{\varphi {\ }|{\ }\varphi \in L_{p}(\mathbb{M}^{d}),{\ }%
\parallel \varphi \parallel _{p}\leq 1\}$ be the unit ball of the space $%
L_{p}(\mathbb{M}^{d})$. The Hilbert space $L_{2}(\mathbb{M}^{d})$ with usual
scalar product $\langle f,g\rangle =\int_{\mathbb{M}^{d}}f(x)\overline{g(x)}%
d\nu (x)$ has the decomposition $L_{2}(\mathbb{M}^{d})=\bigoplus_{k=0}^{%
\infty }\mathrm{H}_{k},$ where $\mathrm{H}_{k}$ is the eigenspace of the
Laplace-Beltrami operator corresponding to the eigenvalue $\theta
_{k}=-k(k+\alpha +\beta +1)$. Let $\{Y_{j}^{k}\}_{j=1}^{d_{k}}$ be an
orthonormal basis of $\mathrm{H}_{k}$. The following addition formula is
known Koornwinder \cite{koor}
\begin{equation}
\sum_{j=1}^{d_{k}}Y_{j}^{k}(x)\overline{Y_{j}^{k}(y)}=\tilde{Z}_{k}(\cos
2\lambda \theta ),  \label{add}
\end{equation}%
where $\theta =d(x,y)$. Comparing (\ref{add}) with (\ref{z}) we get
\begin{equation}
\sum_{j=1}^{d_{k}}Y_{j}^{k}(x)\overline{Y_{j}^{k}(y)}=\tilde{Z}_{k}(\cos
\theta )=C_{k}(\mathbb{M}^{d})P_{k}^{(\alpha ,\beta )}(\cos 2\lambda \theta
).  \label{add1}
\end{equation}


\section{Sets of smooth functions and multiplier operators on $\mathbb{M}%
^{d} $}

Using multiplier operators we introduce a wide range of smooth functions on $%
\mathbb{M}^{d}$. Let $\varphi \in L_{p}(\mathbb{M}^{d})$, $1 \leq p \leq
\infty$, with the formal Fourier expansion
\[
\varphi \sim \sum_{k=0}^{ \infty} \sum_{j=1}^{d_k} c_{k,j}(\phi) Y_{j}^{k},
\,\,\, c_{k,j}(\phi) = \int_{\mathbb{M}^{d}} \phi \overline{Y_{j}^{k}} d\nu.
\]
Let $\Lambda = \{ \lambda_{k} \}_{k \in \mathbb{N}}$ be a sequence of real
(complex) numbers. If for any $\phi \in L_{p}(\mathbb{M}^{d})$ there is a
function $f =\Lambda \phi \in L_{q}(\mathbb{M}^{d})$ such that
\[
f \sim \sum_{k=0}^{ \infty} \lambda_{k}\sum_{j=1}^{d_k} c_{k,j}(\phi)
Y_{j}^{k},
\]
then we shall say that the multiplier operator $\Lambda$ is of $(p, q)$-type
with norm $\| \Lambda \|_{p, q}= \sup_{\varphi \in U_p}\|\Lambda \varphi
\|_{q}$. We shall say that the function $f$ is in $\Lambda U_{p} \oplus
\mathbb{R}$ if
\[
\Lambda \phi = f \sim C {\ }+ \sum_{k=1}^{\infty} \lambda_{k}
\sum_{j=1}^{d_k} c_{k,j}( \phi)Y_{j}^{k},
\]
where $C \in \mathbb{R}$ and $\varphi \in U_p$. In particular, the $\gamma$%
-th fractional integral ($\gamma >0$) of a function $\varphi \in L_{1}(%
\mathbb{M}^{d})$ is defined by the sequence $\lambda_{k}=(k(k + \alpha +
\beta +1))^{- \gamma /2}$. Sobolev's classes $W^{\gamma}_{p}(\mathbb{M}^{d})$
on $\mathbb{M}^{d}$ are defined as sets of functions with formal Fourier
expansions
\[
C + \sum_{k=1}^{ \infty} (k(k + \alpha + \beta +1))^{- \gamma /2}
\sum_{j=1}^{d_k} c_{k,j}( \phi) Y_{j}^{k},
\]
where $C \in \mathbb{R}$ and $\| \phi \|_{p} \leq 1$. Let $Z$ be a zonal
integrable function on $\mathbb{M}^{d}$. For any integrable function $g$ we
can define convolution $h$ on $\mathbb{M}^{d}$ as the following
\[
h(\cdot) = (Z \ast g)(\cdot) = \int_{\mathbb{M}^{d}} Z(\cos(2\lambda
d(\cdot, x)) g(x) d\nu(x).
\]
For the convolution on $\mathbb{M}^{d}$ we have Young's inequality
$\|(z \ast g)\|_{q} \leq \|z\|_{p}\|g\|_{r}, $ 
where $1/q=1/p+1/r-1$ and $1 \leq p,q,r \leq \infty$. It is possible to show
that for any $\gamma > 0$ the function $G_{\gamma} = G_{\gamma, \eta} \sim
\sum_{k=1}^{\infty} (k(k + \alpha + \beta +1))^{- \gamma /2} Z_{k}^{\eta} $
with the pole $\eta$ is integrable on $\mathbb{M}^{d}$ and for any function $%
g \in W^{\gamma}_{p}(\mathbb{M}^{d})$ we have an integral representation $g
= C + G_{\gamma} \ast \phi, $ where $C \in \mathbb{R}$ and $\phi \in U_{p}$.


\section{The Orthogonal Projection}

The main result of this article is the following statement.

\textbf{Theorem 1.} \emph{Let $\mathbb{M}^{d}=\mathbb{S}^{d}$, $\mathrm{P}%
^{d}(\mathbb{C})$, $\mathrm{P}^{d}(\mathbb{H})$, $\mathrm{P}^{16}(\mathrm{Cay%
})$, $d\geq 2$, then
\[
\Vert S_{n}\Vert _{C(\mathbb{M}^{d})\rightarrow C(\mathbb{M}^{d})}=\mathcal{K%
}(\mathbb{M}^{d})n^{(d-1)/2}+O\left\{
\begin{array}{cc}
1, & d=2,3 \\
n^{(d-3)/2}, & d\geq 4%
\end{array}%
\right\} ,
\]%
where
\[
\mathcal{K}(\mathbb{M}^{d})=\frac{4}{\pi ^{3/2}\Gamma (d/2)}%
\,\,\int_{0}^{\pi /2}(\sin \,\,\eta )^{(d-3)/2}\,\,(\cos \,\,\eta )^{\chi (%
\mathbb{M}^{d})}d\eta ,
\]%
and
\[
\chi (\mathbb{M}^{d})=\left\{
\begin{array}{cc}
(d-1)/2, & \mathbb{M}^{d}=\mathbb{S}^{d},\,\,d=2,3,4,\cdots , \\
1/2, & \mathbb{M}^{d}=\mathrm{P}^{d}(\mathbb{C}),\,\,d=4,6,8,\cdots , \\
2, & \mathbb{M}^{d}=\mathrm{P}^{d}(\mathbb{H}),\,\,d=8,12,16,\cdots , \\
7/2, & \mathbb{M}^{d}=\mathrm{P}^{16}(\mathrm{Cay}).%
\end{array}%
\right.
\]%
If $\mathbb{M}^{d}=\mathrm{P}^{d}(\mathbb{R})$, $d=2,3,\cdots ,$ then
\[
\Vert S_{2n}\Vert _{C(\mathrm{P}^{d}(\mathbb{R}))\rightarrow C(\mathrm{P}%
^{d}(\mathbb{R}))}
\]%
\[
=\frac{4\,n^{(d-1)/2}}{\pi ^{3/2}\,\,\Gamma (d/2)}\int_{0}^{\pi /2}(\sin
\eta )^{(d-3)/2}d\eta +O\left\{
\begin{array}{cc}
1, & d=2,3 \\
n^{(d-3)/2}, & d\geq 4.%
\end{array}%
\right\} .
\]%
} 

\textbf{Proof.} Consider the case $\mathbb{M}^{d}=\mathbb{S}^{d},\mathrm{P}%
^{d}(\mathbb{C}),\mathrm{P}^{d}(\mathbb{H}),\mathrm{P}^{16}(\mathrm{Cay})$
first. We will need an explicit representation for the constant $C_{k}(%
\mathbb{M}^{d})$ defined in (\ref{add1}) for our applications. Putting $y=x$
in (\ref{add1}) and then integrating both sides with respect to $d\nu (x)$
we get
\begin{equation}
d_{k}=\dim H_{k}=\sum_{j=1}^{d_{k}}\int_{\mathbb{M}^{d}}|Y_{j}^{k}(x)|^{2}d%
\nu (x)=C_{k}(\mathbb{M}^{d})P_{k}^{(\alpha ,\beta )}(1).  \label{dimhn1}
\end{equation}%
Taking the square of both sides of (\ref{add1}) and then integrating with
respect to $d\nu (x)$ we find
\begin{equation}
\sum_{j=1}^{d_{k}}|Y_{j}^{k}(y)|^{2}=C_{k}^{2}(\mathbb{M}^{d})\int_{\mathbb{M%
}^{d}}\left( P_{k}^{(\alpha ,\beta )}(\cos (2\lambda d(x,y))\right) ^{2}d\nu
(x).  \label{111}
\end{equation}%
Since $d\nu $ is shift invariant then
\[
\int_{\mathbb{M}^{d}}\left( P_{k}^{(\alpha ,\beta )}(\cos (2\lambda
d(x,y)))\right) ^{2}d\nu (x)=c^{-1}\Vert P_{k}^{(\alpha ,\beta )}\Vert
_{2}^{2},
\]%
where the constant $c$ is defined by (\ref{c}) and (see \cite[p.68]{sze})
\[
\Vert P_{k}^{(\alpha ,\beta )}\Vert _{2}^{2}=\int_{-1}^{1}\left(
P_{k}^{(\alpha ,\beta )}(t)\right) ^{2}(1-t)^{\alpha }(1+t)^{\beta }dt
\]%
\[
=\frac{2^{\alpha +\beta +1}}{2k+\alpha +\beta +1}\,\,\frac{\Gamma (k+\alpha
+1)\Gamma (k+\beta +1)}{\Gamma (k+1)\Gamma (k+\alpha +\beta +1)}.
\]%
So that, (\ref{111}) can be written in the form
\[
\sum_{j=1}^{d_{k}}|Y_{j}^{k}(y)|^{2}=c^{-1}C_{k}^{2}(\mathbb{M}^{d})\Vert
P_{k}^{(\alpha ,\beta )}\Vert _{2}^{2}.
\]%
Integrating the last line with respect to $d\nu (y)$ we obtain
\[
d_{k}=c^{-1}C_{k}^{2}(\mathbb{M}^{d})\Vert P_{k}^{(\alpha ,\beta )}\Vert
_{2}^{2}.
\]%
It is sufficient to compare this with (\ref{dimhn1}) to obtain
\begin{equation}
C_{k}(\mathbb{M}^{d})=\frac{cP_{k}^{(\alpha ,\beta )}(1)}{\Vert
P_{k}^{(\alpha ,\beta )}\Vert _{2}^{2}}.  \label{opa}
\end{equation}%
We get now an integral representation for the Fourier sums $S_{n}(\phi ,x)$
of a function $\phi \in L_{1}(\mathbb{M}^{d})$,
\[
S_{n}(\phi ,x)=\sum_{k=0}^{n}\sum_{j=1}^{d_{k}}c_{k,j}(\phi )Y_{j}^{k}(x)
\]%
\[
=\sum_{k=0}^{n}\sum_{j=1}^{d_{k}}\left( \int_{\mathbb{M}^{d}}\phi (y)%
\overline{Y_{j}^{k}(y)}d\mu (y)\right) Y_{j}^{k}(x)
\]%
\[
=\int_{\mathbb{M}^{d}}\sum_{k=0}^{n}\left( \sum_{j=1}^{d_{k}}\overline{%
Y_{j}^{k}(y)}Y_{j}^{k}(x)\right) \phi (y)d\nu (y)
\]%
\[
=\int_{\mathbb{M}^{d}}\sum_{k=0}^{n}Z_{k}^{x}(y)\phi (y)d\nu (y)
\]%
\begin{equation}
=\int_{\mathbb{M}^{d}}K_{n}(x,y)\phi (y)d\nu (y),  \label{laplace}
\end{equation}%
where
\begin{equation}
K_{n}(x,y)=\sum_{k=0}^{n}Z_{k}^{x}(y).  \label{000}
\end{equation}%
By (\ref{z}) and (\ref{opa}) we have
\[
K_{n}(x,y)=c\sum_{k=0}^{n}\,\,\frac{P_{k}^{(\alpha ,\beta )}(1)}{\Vert
P_{k}^{(\alpha ,\beta )}\Vert _{2}^{2}}P_{k}^{(\alpha ,\beta )}(\cos
2\lambda d(x,y)).
\]%
Put
\[
G_{n}(\gamma ,\delta )=\sum_{k=0}^{n}\frac{P_{k}^{(\alpha ,\beta )}(\gamma
)P_{k}^{(\alpha ,\beta )}(\delta )}{\Vert P_{k}^{(\alpha ,\beta )}\Vert
_{2}^{2}},
\]%
then Szeg\"{o} \cite[p.71]{sze},
\[
G_{n}(\gamma ,1)=\sum_{k=0}^{n}\frac{P_{k}^{(\alpha ,\beta )}(\gamma
)P_{k}^{(\alpha ,\beta )}(1)}{\Vert P_{k}^{(\alpha ,\beta )}\Vert _{2}^{2}}
\]%
\begin{equation}
=2^{-\alpha -\beta -1}\frac{\Gamma (n+\alpha +\beta +2)}{\Gamma (\alpha
+1)\Gamma (n+\beta +1)}P_{n}^{(\alpha +1,\beta )}(\gamma ).  \label{g}
\end{equation}%
It means that the kernel function (\ref{000}) in the integral representation
(\ref{laplace}) can be written in the form
\[
K_{n}(x,y)=c2^{-\alpha -\beta -1}\frac{\Gamma (n+\alpha +\beta +2)}{\Gamma
(\alpha +1)\Gamma (n+\beta +1)}P_{n}^{(\alpha +1,\beta )}(\cos 2\lambda
d(x,y)).
\]%
Let $\mathbf{o}$ be the north pole of $\mathbb{M}^{d}$, then since $K_{n}$
is a zonal function and $d\nu $ is shift invariant,
\[
\| S_{n}\| _{C(\mathbb{M}^{d})\rightarrow C(\mathbb{M}%
^{d})}=\sup_{\| \phi \| _{C(\mathbb{M}^{d})}\leq 1}\Vert S_{n}(\phi
,x)\| _{C(\mathbb{M}^{d})}
\]%
\[
=\sup \left\{ \int_{\mathbb{M}^{d}}|K_{n}(x,y)|d\nu (y)\,|\,\,x\in \mathbb{M}%
^{d}\right\}
\]%
\[
=\int_{\mathbb{M}^{d}}|K_{n}(\mathbf{o},y)|d\nu (y)
\]%
\[
=c2^{-\alpha -\beta -1}\frac{\Gamma (n+\alpha +\beta +2)}{\Gamma (\alpha
+1)\Gamma (n+\beta +1)}\int_{\mathbb{M}^{d}}|P_{n}^{(\alpha +1,\beta )}(\cos
2\lambda d(\mathbf{o},y))|d\nu (y)
\]%
\[
=cc^{-1}2^{-\alpha -\beta -1}\frac{\Gamma (n+\alpha +\beta +2)}{\Gamma
(\alpha +1)\Gamma (n+\beta +1)}\int_{-1}^{1}|P_{n}^{(\alpha +1,\beta
)}(t)|(1-t)^{\alpha }(1+t)^{\beta }dt
\]%
\[
=\frac{2^{-\alpha -\beta -1}\Gamma (n+\alpha +\beta +2)}{\Gamma (\alpha
+1)\Gamma (n+\beta +1)}\int_{0}^{\pi }|P_{n}^{(\alpha +1,\beta )}(\cos \eta
)|\left( 2\sin ^{2}\frac{\eta }{2}\right) ^{\alpha }\left( 2\cos ^{2}\frac{%
\eta }{2}\right) ^{\beta }\sin \eta d\eta
\]%
\[
=\frac{\Gamma (n+\alpha +\beta +2)}{\Gamma (\alpha +1)\Gamma (n+\beta +1)}%
\int_{0}^{\pi }|P_{n}^{(\alpha +1,\beta )}(\cos \eta )|\left( \sin \frac{%
\eta }{2}\right) ^{2\alpha +1}\left( \cos \frac{\eta }{2}\right) ^{2\beta
+1}d\eta
\]%
\begin{equation}
=\left( \frac{I_{n}}{\Gamma (\alpha +1)}\right) (n^{\alpha +1}+O(n^{\alpha
})),\text{ }n\rightarrow \infty ,  \label{final}
\end{equation}%
where
\begin{equation}
I_{n}=\int_{0}^{\pi }|P_{n}^{(\alpha +1,\beta )}(\cos \eta )|\left( \sin
\frac{\eta }{2}\right) ^{2\alpha +1}\left( \cos \frac{\eta }{2}\right)
^{2\beta +1}d\eta .  \label{i}
\end{equation}%
It is known Szeg\"{o} \cite[p.196]{sze} that for $0<\theta <\pi $,
\begin{equation}
P_{n}^{(\alpha +1,\beta )}(\cos \eta )=n^{-1/2}\,\,\kappa (\eta )\cos (N\eta
+\gamma )+O(n^{-3/2}),  \label{as}
\end{equation}%
where
\[
\kappa (\eta )=\pi ^{-1/2}\left( \sin \frac{\eta }{2}\right) ^{-\alpha
-3/2}\left( \cos \frac{\eta }{2}\right) ^{-\beta -1/2}
\]%
and
\[
N=n+1+\frac{\alpha +\beta }{2},\,\,\,\gamma =-\frac{\alpha +3/2}{2}\pi =-%
\frac{d+1}{4}\pi .
\]%
Comparing (\ref{final}) - (\ref{as}), applying a simple Tylor series
arguments and elementary estimates of the derivative of the function
\[
\sigma (\eta )\,=\,\left( \sin \frac{\eta }{2}\right) ^{\alpha -1/2}\left(
\cos \frac{\eta }{2}\right) ^{\beta +1/2}
\]%
we get
\[
I_{n}\,=\,\pi ^{-1/2}n^{-1/2}\int_{0}^{\pi }\left( \sin \frac{\eta }{2}%
\right) ^{\alpha -1/2}\left( \cos \frac{\eta }{2}\right) ^{\beta +1/2}
\]%
\[
\times \left\| \,\,\cos \left( \left( n+\frac{\alpha +\beta +2}{2}\right)
\eta -\frac{d+1}{4}\pi \right) \right\| d\eta +\,\,O(n^{-3/2})
\]%
\[
=2\pi ^{-3/2}n^{-1/2}\int_{0}^{\pi }\left( \sin \frac{\eta }{2}\right)
^{\alpha -1/2}\left( \cos \frac{\eta }{2}\right) ^{\beta +1/2}d\eta
\]%
\begin{equation}
+n^{-1/2}\,O\,\left\{
\begin{array}{cc}
n^{-1/2}, & \alpha =0 \\
n^{-1}, & \alpha \geq 1/2%
\end{array}%
\right\} ,\text{ \ }n\rightarrow \infty .  \label{sao}
\end{equation}%
Remind that $\alpha =(d-2)/2$, $d\geq 2$ for any manifold $\mathbb{M}^{d}$
under consideration. Put $\chi (\mathbb{M}^{d})=\beta +1/2$, then from (\ref%
{final}) and (\ref{sao}) it follows that
\[
\| S_{n}\| _{C(\mathbb{M}^{d})\rightarrow C(\mathbb{M}^{d})}=\mathcal{K%
}(\mathbb{M}^{d}n^{\alpha +1/2}+O(n^{\alpha -1/2}),
\]%
where
\[
\mathcal{K}(\mathbb{M}^{d})=\frac{2}{\pi ^{3/2}\Gamma (\alpha +1)}%
\,\,\int_{0}^{\pi }\left( \sin \frac{\eta }{2}\right) ^{\alpha -1/2}\left(
\cos \frac{\eta }{2}\right) ^{\beta +1/2}d\eta
\]%
%
%
%
%
%
%
%
%
%
\[
=\frac{4}{\pi ^{3/2}\Gamma (d/2)}\,\,\int_{0}^{\pi /2}(\sin \eta
)^{(d-3)/2}(\cos \eta )^{\chi (\mathbb{M}^{d})}d\eta ,
\]%
since $\alpha =(d-2)/2$. Hence,
\[
\| S_{n}\| _{C(\mathbb{M}^{d})\rightarrow C(\mathbb{M}^{d})}=\mathcal{K%
}(\mathbb{M}^{d})n^{(d-1)/2}+O\,\left\{
\begin{array}{cc}
n^{(d-2)/2}, & d=2 \\
n^{(d-3)/2}, & d\geq 3%
\end{array}%
\right\} .
\]%
Finally, the value of $\chi (\mathbb{M}^{d})$, where $\mathbb{M}^{d}=\mathbb{%
S}^{d},\mathrm{P}^{d}(\mathbb{C}),\mathrm{P}^{d}(\mathbb{H}),\mathrm{P}^{16}(%
\mathrm{Cay})$, can be easily calculated using (\ref{pipa}) and (\ref{pipo}).

The case of $\mathrm{P}^{d}(\mathbb{R})$ needs a special treatment. In this
case $\alpha =\beta =(d-2)/2$, $\lambda =\pi /(4L)$ and the kernel function $%
K_{2n}^{\ast }(x,y)$ in the integral representation for the Fourier sums,
\[
S_{2n}(\phi ,x)=\int_{\mathrm{P}^{d}(\mathbb{R})}K_{2n}^{\ast }(x,y)\phi
(y)d\nu (y)
\]%
has the form
\[
K_{2n}^{\ast }(x,y)=\sum_{k=0}^{n}Z_{2k}^{x}(y)=\sum_{k=0}^{n}C_{2k}(\mathrm{%
P}^{d}(\mathbb{R}))P_{2k}^{(\alpha ,\alpha )}\cos (2\lambda d(x,y))
\]%
\[
=\sum_{k=0}^{n}C_{2k}(\mathrm{P}^{d}(\mathbb{R}))P_{2k}^{((d-2)/2,(d-2)/2)}%
\left( \cos \left( \frac{\pi }{2L}d(x,y)\right) \right) .
\]%
Let the constant $c^{\ast }$ be such that
\[
1=\int_{\mathrm{P}^{d}(\mathbb{R})}d\nu =\int_{0}^{1}\omega
^{(d-2)/2,(d-2)/2}(t)dt=(c^{\ast })^{-1}\int_{0}^{1}(1-t^{2})^{(d-2)/2}dt,
\]%
then $c^{\ast }=c/2$ and
\[
C_{2k}(\mathrm{P}^{d}(\mathbb{R}))\,=\,\frac{c^{\ast
}P_{2k}^{((d-2)/2,(d-2)/2)}(1)}{\left\|
P_{2k}^{((d-2)/2,(d-2)/2)}\right\| _{2,\ast }^{2}}=\frac{%
cP_{2k}^{((d-2)/2,(d-2)/2)}(1)}{\left\|
P_{2k}^{((d-2)/2,(d-2)/2)}\right\| _{2}^{2}},
\]%
where
\[
\left\| P_{2k}^{((d-2)/2,(d-2)/2)}\right\| _{2\ast
}^{2}=\int_{0}^{1}\left( P_{2k}^{((d-2)/2,(d-2)/2)}(t)\right)
^{2}(1-t^{2})^{(d-2)/2}
\]%
\[
=2^{-1}\,\| P_{2k}^{((d-2)/2,(d-2)/2)}\| _{2}^{2}.
\]%
Let $\mathbf{o}$ be the north pole of $\mathrm{P}^{d}(\mathbb{R})$, then
since $K_{2n}^{\ast }$ is a zonal function and $d\nu $ is shift invariant,
\[
\| S_{2n}\| _{C(\mathrm{P}^{d}(\mathbb{R}))\rightarrow C(\mathrm{P}%
^{d}(\mathbb{R}))}=\sup_{\| \phi \| _{C(\mathrm{P}^{d}(\mathbb{R}%
))}\leq 1}\| S_{2n}(\phi )\| _{C(\mathrm{P}^{d}(\mathbb{R}))}
\]%
\begin{equation}
=\sup \left\{ \int_{\mathrm{P}^{d}(\mathbb{R})}|K_{2n}^{\ast }(x,y)|d\nu
(y)|\,\,x\in \mathrm{P}^{d}(\mathbb{R})\right\} =\int_{\mathrm{P}^{d}(%
\mathbb{R})}|K_{2n}^{\ast }(\mathbf{o},y)|d\nu (y).  \label{kkk}
\end{equation}%
Consider the function
\[
G_{2n}^{\ast }(\gamma ,1)=\sum_{k=0}^{n}\frac{P_{2k}^{(\alpha ,\alpha
)}(\gamma )P_{2k}^{(\alpha ,\alpha )}(1)}{\| P_{2k}^{(\alpha ,\alpha
)}\| _{2}^{2}}.
\]%
Since $P_{k}^{(\alpha ,\beta )}(\gamma )=(-1)^{k}P_{k}^{(\beta ,\alpha
)}(-\gamma )$, Szeg\"{o} \cite[p.59]{sze}, then
\[
G_{2n}^{\ast }(\gamma ,1)=\frac{G_{2n}(\gamma ,1)+G_{2n}(-\gamma ,1)}{2}
\]%
\[
=\frac{2^{-d+1}\Gamma (2n+d)}{\Gamma (d/2)\Gamma (2n+d/2)}\left( \frac{%
P_{2n}^{(d/2,(d-2)/2)}(\gamma )+P_{2n}^{(d/2,(d-2)/2)}(-\gamma )}{2}\right)
\]%
\[
=\frac{2^{-d}\Gamma (2n+d)}{\Gamma (d/2)\Gamma (2n+d/2)}\left(
P_{2n}^{(d/2,(d-2)/2)}(\gamma )+P_{2n}^{((d-2)/2,d/2)}(\gamma )\right)
\]%
where $G_{2n}(\gamma ,1)$ is defined in (\ref{g}). Consequently, (\ref{kkk})
takes the form
\[
\| S_{2n}\| _{C(\mathrm{P}^{d}(\mathbb{R}))\rightarrow C(\mathrm{P}%
^{d}(\mathbb{R}))}
\]%
\[
=\frac{c\,2^{-d}\Gamma (2n+d)}{\Gamma (d/2)\Gamma (2n+d/2)}
\]%
\[
\times \int_{\mathrm{P}^{d}(\mathbb{R})}\left\|
P_{2n}^{(d/2,(d-2)/2)}(\cos (\pi d(\mathbf{o},y)/(4L)))\,+%
\,P_{2n}^{((d-2)/2,d/2)}(\cos (\pi d(\mathbf{o},y)/(4L)))\right\| d\nu (y)
\]%
\[
=\frac{c\,2^{-d}\Gamma (2n+d)}{c^{\ast }\Gamma (d/2)\Gamma (2n+d/2)}%
\int_{0}^{1}\left\|
P_{2n}^{(d/2,(d-2)/2)}(t)+P_{2n}^{((d-2)/2,d/2)}(t)\right\|
(1-t^{2})^{(d-2)/2}dt
\]%
\begin{equation}
=\frac{2^{-d+1}\Gamma (2n+d)}{\Gamma (d/2)\Gamma (2n+d/2)}I_{n}^{^{\prime }},
\label{in}
\end{equation}%
where
\[
I_{n}^{^{\prime }}=\int_{0}^{1}\left\|
P_{2n}^{(d/2,(d-2)/2)}(t)+P_{2n}^{((d-2)/2,d/2)}(t)\right\|
(1-t^{2})^{(d-2)/2}dt
\]%
\[
=\int_{0}^{\pi /2}\left\| P_{2n}^{(d/2,(d-2)/2)}(\cos \eta
)+P_{2n}^{((d-2)/2,d/2)}(\cos \eta )\right\| (\sin \eta )^{d-1}d\eta
\]%
Applying (\ref{as}) we get
\[
I_{n}^{^{\prime }}=\frac{1}{\pi ^{1/2}2^{1/2}n^{1/2}}\int_{0}^{\pi /2}d\eta
(\sin \eta )^{d-1}
\]%
\[
\times \left\| \left( \sin \frac{\eta }{2}\right) ^{-d/2-1/2}\right.
\left( \cos \frac{\eta }{2}\right) ^{-(d-2)/2-1/2}
\]%
\[
\times \cos \left. \left( \left( 2n+\frac{d/2+(d-2)/2+1}{2}\right) \eta -%
\frac{(d/2+1/2)\pi }{2}\right) \right.
\]%
\[
+\left( \sin \frac{\eta }{2}\right) ^{-(d-2)/2-1/2}\left( \cos \frac{\eta }{2%
}\right) ^{-d/2-1/2}
\]%
\[
\times \cos \left. \left( \left( 2n+\frac{(d-2)/2+d/2+1}{2}\right) \eta -%
\frac{((d-2)/2+1/2)\pi }{2}\right) \right\| +O(n^{-3/2})
\]%
\[
=\frac{1}{\pi ^{1/2}2^{1/2}n^{1/2}}\int_{0}^{\pi /2}(\sin \eta )^{d-1}\left(
\sin \frac{\eta }{2}\right) ^{-d/2-1/2}\left( \cos \frac{\eta }{2}\right)
^{-d/2-1/2}
\]%
\[
\times \left\| \cos \frac{\eta }{2}\,\cos \left( \left( 2n+\frac{d}{2}%
\right) \eta -\frac{(d+1)\pi }{4}\right) +\sin \frac{\eta }{2}\,\cos \left(
\left( n+\frac{d}{2}\right) \eta -\frac{(d-1)\pi }{4}\right) \right\|
d\eta
\]%
\[
+O(n^{-3/2})
\]%
\[
=\frac{2^{d/2+1/2}}{\pi ^{1/2}2^{1/2}n^{1/2}}\int_{0}^{\pi /2}(\sin \eta
)^{(d-3)/2}\left\| \cos \left( \left( 2n+\frac{d-1}{2}\right) \eta -\frac{%
(d+1)\pi }{4}\right) \right\| d\eta +O(n^{-3/2})
\]%
\[
=\frac{2^{d/2+1}}{\pi ^{3/2}n^{1/2}}\int_{0}^{\pi /2}(\sin \eta
)^{(d-3)/2}d\eta +O\left\{
\begin{array}{cc}
n^{-1/2}, & d=2 \\
n^{-1}, & d\geq 3%
\end{array}%
\right\} .
\]%
Comparing (\ref{in}) with the last line we get
\[
\| S_{2n}\| _{C(\mathrm{P}^{d}(\mathbb{R}))\rightarrow C(\mathrm{P}%
^{d}(\mathbb{R}))}=\mathcal{K}(\mathrm{P}^{d}(\mathbb{R}))\,\,n^{(d-1)/2}+O%
\left\{
\begin{array}{cc}
n^{(d-2)/2}, & d=2 \\
n^{(d-3)/2}, & d\geq 3%
\end{array}%
\right\} ,
\]%
where
\[
\mathcal{K}(\mathrm{P}^{d}(\mathbb{R}))=\frac{4}{\pi ^{3/2}\,\,\Gamma (d/2)}%
\int_{0}^{\pi /2}(\sin \eta )^{(d-3)/2}d\eta .
\]%
$\square $

\textbf{Remark 1.} \emph{Let $\mathbb{M}^{d}=\mathbb{S}^{d},\mathrm{P}^{d}(%
\mathbb{R}),\mathrm{P}^{d}(\mathbb{C}),\mathrm{P}^{d}(\mathbb{H}),\mathrm{P}%
^{16}(\mathrm{Cay})$. It is known \cite{jfan} that for any $\gamma >0$,
\[
E_{n}(W_{\infty }^{\gamma }(\mathbb{M}^{d}))=\sup \{E_{n}(f)|\,\,f\in
W_{\infty }^{\gamma }(\mathbb{M}^{d})\}\asymp n^{-\gamma }.
\]%
From the Theorem 1 and the Lebesgue inequality it follows that the Fourier
series of a function $f\in W_{\infty }^{\gamma }(\mathbb{M}^{d})$ converges
uniformly if $\gamma >(d-1)/2$. In general, let $\Delta ^{0}\lambda
_{k}=\lambda _{k}$, $\Delta ^{1}\lambda _{k}=\lambda _{k}-\lambda _{k+1}$, $%
\Delta ^{s+1}\lambda _{k}=\Delta ^{s}\lambda _{k}-\Delta ^{s}\lambda _{k+1}$%
, \thinspace \thinspace $k,s\in \mathbb{N}$ and
\[
M:=\left\{
\begin{array}{cc}
(d+1)/2, & d=3,5,\cdots , \\
(d+2)/2, & d=2,4,\cdots%
\end{array}%
\right.
\]%
Let $\Lambda =\{\lambda _{k}\}_{k\in \mathbb{N}}$ be a multiplier operator, $%
\Lambda :L_{\infty }(\mathbb{M}^{d})\rightarrow L_{\infty }(\mathbb{M}^{d})$
and $\Lambda U_{\infty }(\mathbb{M}^{d})$ be the respective set of smooth
functions, then from the Theorem 2, \cite[p.317]{jfan} it follows that the
Fourier series of a function $f\in \Lambda U_{\infty }(\mathbb{M}^{d})$
converges uniformly if
\[
\lim_{n\rightarrow \infty }n^{(d-1)/2}\sum_{k=n+1}^{\infty }|\Delta
^{M+1}\lambda _{k}|\,\,k^{M}=0,
\]%
since $E_{n}(\Lambda U_{\infty }(\mathbb{M}^{d}))\ll \sum_{k=n+1}^{\infty
}|\Delta ^{M+1}\lambda _{k}|\,\,k^{M}$. In particular, let
\[
\Lambda \,=\,\{\lambda _{k}\}_{k\in \mathbb{N}},\,\,\,\lambda
_{k}\,=\,k^{-(d-1)/2}\,(\mathrm{ln}\,k)^{-\alpha },
\]%
where $\alpha >0$, then the Fourier series of any function $f\in \Lambda
U_{\infty }(\mathbb{M}^{d})$ converges uniformly. A similar result is valid
for $\mathrm{P}^{d}(\mathbb{R})$. }

\textbf{Remark 2.} \emph{In terms of gamma function we have the following
representations for the constant $\mathcal{K}(\mathbb{M}^{d})$:
\[
\mathcal{K}(\mathbb{S}^{d})= \frac{2\,\Gamma\left(\frac{d-1}{4}\right)
\Gamma\left(\frac{d+1}{4}\right)}{\pi^{3/2} \left(\Gamma\left(\frac{d}{2}%
\right)\right)^{2}},\,\,\,d=2,3,4,\cdots
\]
\[
\mathcal{K}(\mathrm{P}^{d}(\mathbb{R}))= \frac{2\,\Gamma\left(\frac{d-1}{4}%
\right)}{\pi \Gamma\left(\frac{d}{2}\right)\Gamma\left(\frac{d+1}{4}\right)}%
,\,\,\,d=2,3,4,\cdots
\]
\[
\mathcal{K}(\mathrm{P}^{d}(\mathbb{C}))= \frac{2\,\Gamma\left(\frac{d-1}{4}%
\right) \Gamma\left(\frac{3}{4}\right)}{\pi^{3/2} \Gamma\left(\frac{d}{2}%
\right)\Gamma\left(\frac{d+2}{4}\right)},\,\,\,d=4,6,8,\cdots
\]
\[
\mathcal{K}(\mathrm{P}^{d}(\mathbb{H}))= \frac{\Gamma\left(\frac{d-1}{4}%
\right)}{\pi \Gamma\left(\frac{d}{2}\right)\Gamma\left(\frac{d+5}{4}\right)}%
,\,\,\,d=8,12,16,\cdots
\]
\[
\mathcal{K}(\mathrm{P}^{16}(\mathrm{Cay}))= \frac{11 \cdot 2^{1/2}}{2949120
\cdot \pi^{1/2}}.
\]
}

Here we present the solution of the problem of Kolmogorov on sharp
asymptotic of the rate of convergence of Fourier series on sets of smooth
functions on manifolds. First we will need some definitions.

For a given multiplier sequence $\Lambda =\{\lambda _{k}\}$ let
\[
C_{n}^{\delta }=\frac{\Gamma (n+\delta +1)}{\Gamma (\delta +1)\Gamma (n+1)}%
\asymp n^{\delta }
\]%
and
\[
S_{n}^{\delta }=\frac{1}{C_{n}^{\delta }}\sum_{m=0}^{n}C_{n-m}^{\delta
}Z_{m}.
\]%
It is known \cite{jfan} that
\begin{equation}
\| S_{n}^{\delta }\| _{1}\ll \left\{
\begin{array}{cc}
1, & \delta >(d-1)/2, \\
\log \,n, & \delta =(d-1)/2, \\
n^{(d-1)/2-\delta }, & 0\leq \delta <(d-1)/2.%
\end{array}%
\right.  \label{delta}
\end{equation}

\textbf{Definition 1.} \emph{Let $\mathcal{T}$ be the set of multipliers $%
\Lambda =\{\lambda _{k}\}_{k\in \mathbb{N}}$ such that
\[
\lim_{m\rightarrow \infty }\left\|\sum_{s=0}^{d-1}\Delta ^{s}\lambda
_{m-s}C_{m-s}^{s}S_{m-s}^{s}\right\| _{1}=0,
\]%
and
\[
\left\| \sum_{k=n+1}^{\infty }\Delta ^{d}\lambda
_{k}C_{k}^{d-1}S_{k}^{d-1}+\sum_{s=1}^{d}\Delta ^{s}\lambda
_{n+1}C_{n}^{s}S_{n}^{s}\right\| _{1}=o(|\lambda _{n+1}|\|
S_{n}^{0}\| _{1}).
\]%
} 
As a simple consequence of Theorem 1 we get

\textbf{Theorem 2.} \emph{Let $\Lambda \in \mathcal{T}$ then}%
\begin{equation}  \label{101111}
\sup_{f\in \Lambda U_{p}(\mathbb{M}^{d})}\,\| f-S_{n}(f)\| _{p}=%
\mathcal{K}(\mathbb{M}^{d})|\lambda _{n+1}|n^{(d-1)/2}(1+o(1)),
\end{equation}%
where $p=1,\infty ,\,\,\,\Lambda \in \mathcal{T},\,\,\,\mathbb{M}^{d}=%
\mathbb{S}^{d},\,\mathrm{P}^{d}(\mathbb{C}),\,\mathrm{P}^{d}(\mathbb{H}),\,%
\mathrm{P}^{d}(\mathrm{Cay})$ and
\begin{equation}  \label{02}
\sup_{f\in \Lambda U_{p}(\mathrm{P}^{d}(\mathbb{R}))}\,\|
f-S_{2n}(f)\| _{p}=\mathcal{K}(\mathrm{P}^{d}(\mathbb{R}))|\lambda
_{n+1}|n^{(d-1)/2}(1+o(1)),\,\,\,p=1,\infty .
\end{equation}

\textbf{Proof} Since $\Lambda \in \mathcal{T}$ then
\[
\left\|\sum_{k=n+1}^{\infty} \lambda_{k} Z_{k} \right\|_{1}= |\lambda_{n+1}|
\|S_{n}^{0}\|_{1}(1+o(1)).
\]
If $\mathbb{M}^{d}=\mathrm{P}^{d}(\mathbb{R})$ and $\{\lambda_{k}\}_{k \in
\mathbb{N}} \in \mathcal{T}$ then
\[
\left\|\sum_{k=n+1}^{\infty} \lambda_{k} Z_{2k}\right\|_{1} =
|\lambda_{n+1}|\|S^{0}_{n}\|_{1}(1+o(1)),
\]
where
\[
S^{\delta}_{n} = \frac{1}{C^{\delta}_{n}}\,\sum_{m=0}^{n}\,C^{\delta}_{n-m}%
\,Z_{2m},\,\,\,\delta>0.
\]
Hence, applying Theorem 1 we get (23) and (\ref{02}).

\textbf{Remark 3.} \emph{In particular, let $\lambda _{k}=(k(k+\alpha +\beta
+1))^{-\gamma /2}$, $\gamma >0$, then using (\ref{delta}) it is easy to show
that for any fixed $\gamma >0$, $\Lambda =\{\lambda _{k}\}_{n\in \mathbb{N}%
}\in \mathcal{T}$ and
\[
\sup_{f\in W_{p}^{\gamma }(\mathbb{M}^{d})}\,\| f-S_{n}(f)\| _{p}
\]%
}%
\[
=\mathcal{K}(\mathbb{M}^{d}) n^{-\gamma +(d-1)/2}+
O\left( n^{-\gamma }\,\left\{
\begin{array}{cc}
1, & d=2 \\
\ln n, & d=3 \\
n^{(d-3)/2}, & d\geq 4%
\end{array}%
\right\} \right) ,
\]%
\emph{%
\[
p=1,\infty ,\,\,\,\Lambda \in \mathcal{T},\,\,\,\mathbb{M}^{d}=\mathbb{S}%
^{d},\,\mathrm{P}^{d}(\mathbb{C}),\,\mathrm{P}^{d}(\mathbb{H}),\,\mathrm{P}%
^{d}(\mathrm{Cay}),
\]%
\[
\sup_{f\in W_{p}^{\gamma }(\mathrm{P}^{d}(\mathbb{R}))}\,\|
f-S_{2n}(f)\| _{p}
\]%
}%
\[
=\mathcal{K}(\mathrm{P}^{d}(\mathbb{R}))\,n^{-\gamma +(d-1)/2}+O\left(
n^{-\gamma }\,\left\{
\begin{array}{cc}
1, & d=2 \\
\ln n, & d=3 \\
n^{(d-3)/2}, & d\geq 4%
\end{array}%
\right\} \right) .
\]
The case $\mathbb{M}^{d}=\mathbb{S}^{d}$ was considered by C. Xirong, D. Feng, W. Kunyang (Estimations of the remainder
of spherical harmonic series, Math. Proc. Cambridge Philos. Soc. 145 (2008),
no. 1, 243-255).
Namely, it was shown that for any $\gamma >0$ and $d\geq 2$,
\[
\sup_{f\in W_{\infty }^{\gamma }(\mathbb{S}^{d})}\,\,\| f-S_{n}(f)\|
_{\infty }=K_{d}\,n^{-\gamma +\frac{d-1}{2}}+H(n,\gamma ,d)
\]%
where
\[
K_{d}=\frac{2^{\frac{d+1}{2}+2}}{\pi ^{3/2}(d-1)}\,\frac{\left( \Gamma
\left( \frac{d+1}{2}\right) \right) ^{2}}{(d-1)!\Gamma \left( \frac{d}{2}%
\right) }
\]%
\[
=\frac{8\Gamma \left( \frac{d+1}{2}\right) }{2^{(d+1)/2}\pi \Gamma \left(
\frac{d+2}{2}\right) \Gamma \left( \frac{d}{2}\right) }=\frac{8}{%
2^{(d+1)/2}\Gamma \left( \frac{d+2}{2}\right) }\,\frac{\omega _{d}}{\omega
_{d+1}},
\]%
$\omega _{d+1}=\int_{\mathbb{S}^{d}}\,d\sigma (x)$, and $H(n,\gamma ,d)$
satisfies the condition
\[
|H(n,\gamma ,d)|\leq \left\{
\begin{array}{cc}
Cn^{-\gamma }, & if\,\,d=2; \\
Cn^{-\gamma }\log n, & if\,\,d=3; \\
Cn^{-\gamma +\frac{d-3}{2}}, & if\,\,d\geq 4.%
\end{array}%
\right.
\]%
Unfortunately this result is incorrect. In particular, for any $d\geq 2$,
\[
\frac{2^{\frac{d+1}{2}+2}}{\pi ^{3/2}(d-1)}\,\frac{\left( \Gamma \left(
\frac{d+1}{2}\right) \right) ^{2}}{(d-1)!\Gamma \left( \frac{d}{2}\right) }%
\neq \frac{8\Gamma \left( \frac{d+1}{2}\right) }{2^{(d+1)/2}\pi \Gamma
\left( \frac{d+2}{2}\right) \Gamma \left( \frac{d}{2}\right) }.
\]

\vspace{1cm}




\begin{thebibliography}{20}
\bibitem{jfan} B. Bordin, A. Kushpel, J. Levesley, S. Tozoni, Estimates of $%
n $-Widths of Sobolev's Classes on Compact Globally Symmetric Spaces of Rank
1, J. Funct. Anal. 202 (2003) 307-326.


\bibitem{car} E. Cartan, Sur la determination d'un systeme orthogonal
complet dans un espace de Riemann symetrique clos, Rendiconti Circ. mat. di
Palermo, 53 (1929) 217-252.

\bibitem{erd} A. Erd\'elyi, (eds.), Higher Transcendental Functions, Vol. 2,
McGraw-Hill, New York, 1953.

\bibitem{fejer} L. Fejer, Lebesguesche Konstanten und divergente
Fourierreihen, J. f\"ur reine und angew. Math. 138 (1910) 22-53.

\bibitem{gang} R. Gangolli, Positive definite kernels on homogeneous spaces
and certain stochastic processes related to L\'{e}vy's Browian motion of
several parameters, Ann. Inst. H. Poincar\'{e} 3 (1967)(2) 121--225.

\bibitem{gronwall} T.H. Gronwall, On the degree of convergence of Laplace
series, Trans. Amer. Math. Soc. 15 (1914)(1) 1-30.

\bibitem{hel1} S. Helgason, The Radon Transform on Euclidean spaces, compact
two-point homogeneous spaces and Grassmann manifolds, Acta mat. 113 (1965)
153-180.

\bibitem{helgason} S. Helgason, (eds.), Differential geometry and symmetric
spaces", Academic Press, New York, 1962.

\bibitem{koor} T. Koornwinder, The addition formula for Jacobi Polynomials
and spherical harmonics, SIAM J. Appl. Math. 25 (1973)(2) 236-246.

\bibitem{lebesgue} H. Lebesgue, Sur les int\'egrales singuli\`eres, Ann. de
Toulouse 1 (1909) 25-117.

\bibitem{sze} G. Szeg\"{o}(eds.), Orthogonal Polynomials, AMS, New York,
1939.

\bibitem{wang} H.C. Wang, Two-point homogeneous spaces, Annals of Math. 55
(1952) 177-191.


\end{thebibliography}
\end{document}